\documentstyle[fleqn,12pt]{article}
\begin{document}
\pagenumbering{arabic}
\setcounter{page}{1}
\pagestyle{plain}
\baselineskip=16pt

\thispagestyle{empty}
\rightline{MSUMB 96-04, November 1996} 
\vspace{2cm}

\begin{center}
{\bf Two-parametric extension of $h$-deformation of $GL(1|1)$ }
\end{center}

\vspace{1cm}
\noindent
Salih Celik \footnote{New E-mail address: sacelik@yildiz.edu.tr}

\noindent
Mimar Sinan University, Department of Mathematics, 
80690 Besiktas, Istanbul, TURKEY.

\noindent
E-mail: celik@msu.edu.tr 

\vspace{2.5cm}
\noindent
{\bf Abstract}. 
The two-parametric quantum deformation of the algebra of coordinate 
functions on the supergroup GL$(1|1)$ via a contraction 
of GL$_{p,q}(1|1)$ is presented. Related differential calculus on the 
quantum superplane is introduced. 

\noindent
{\bf Mathematics Subject Classifications (1991):} 16S80, 81R50. 

\noindent
{\bf Key words:} supergroup, $q$-deformation, $h$-deformation, differential 
calculus.

\vfill\eject
\noindent
{\bf 1. Introduction }

\noindent
Recently matrix groups like GL$(2)$, GL$(1|1)$, {\it etc.}, 
were generalized in two ways called the $q$-deformation and the 
$h$-deformation. One of them [$q$-deformation] is based on a 
deformation of the algebra of functions on the groups generated 
by coordinate functions $t^i{}_j$ which normally commute. 

In the $q$-deformation of matrix groups [1,2], the commutation 
relations satisfied between the coordinate functions on the 
groups are determined by a matrix $R_q$ so that the functions 
do not commute but satisfy the equation 
$$R_q (T \otimes T) = (I \otimes T)(T \otimes I) R_q,$$
or an endomorphism of the space with non-commuting coordinates 
may be used [3]. In this equation the elements of $R_q$ are 
complex numbers but the matrix $T = (t^i{}_j)$ is formed by 
generally non-commuting elements of an algebra. 

The second type of deformation, called the $h$-deformation, 
is a new class of quantum deformations of matrix groups. 
The $h$-deformation is obtained as a contraction from the 
$q$-deformation of matrix groups and it has been intensively 
studied by many authors [4-12]. 

The purpose of this paper is to present the two-parametric 
extension of the $h$-deformation of the simplest supergroup 
GL$(1|1)$. The single parameter $h$-deformation of GL$(1|1)$, 
GL$_h(1|1)$, was introduced by Dabrowski and Parashar [13]. 
An interesting feature of the two-parameter deformation is 
that both the deformation parameters are anticommuting 
grassmann numbers. 

The paper is organized as follows. In section 2 we give some 
notations and useful formulas which will be used in this work. 
In the following section we present the two-parameter deformation 
of GL$(1|1)$ as related to superplanes. A two parameter R-matrix 
which deforms the supergroup GL$(1|1)$ is introduced in section 4. 
In section 5 we construct the differential calculus on the quantum 
superplane on which the two-parameter quantum supergroup acts. 

%\vfill\eject
\noindent
{\bf 2. Review of GL$_{p,q}(1|1)$}

\noindent
We know that the supergroup GL$(1|1)$ can be deformed by assuming that the 
linear transformations in GL$(1|1)$ are invariant under the action of the 
quantum superplane and its dual [3]. 

In this paper we denote $(p,q)$-deformed objects by primed quantities. 
Unprimed quantities represent transformed coordinates. 

Consider the Manin quantum superplane $A_p$ and its dual $A^*_q$. 
The quantum superplane $A_p$ is generated by coordinates $x'$ and $\xi'$, 
and the commutation rules 
$$ x' \xi' - p \xi' x' = 0, \quad \xi'^2 = 0. \eqno(1) $$
The quantum (dual) superplane $A^*_q$ is generated by coordinates $\eta'$ and 
$y'$, and the commutation rules 
$$ \eta'^2 = 0, \quad \eta' y' - q^{-1} y' \eta' = 0.    \eqno(2) $$
Taking 
$$T' = \left(\matrix{ a' & \beta' \cr \gamma' & d' \cr} \right)$$
as a supermatrix in GL$(1|1)$, we demand that the relations (1), (2) are 
preserved under the action of $T'$ on the quantum superplane and its dual. 
Under some assumption one obtains the following $(p,q)$-commutation 
relations [3,14,15] 
$$ a' \beta' = q \beta' a', \quad a' \gamma' = p\gamma' a', 
   \quad \beta'^2 = 0, $$
$$ d' \beta' = q\beta' d', \quad d' \gamma' = p \gamma' d', 
   \quad \gamma'^2 = 0, \eqno(3)$$
$$\beta' \gamma' + pq^{-1} \gamma' \beta' = 0, \quad 
   a'd' = d'a' + (p - q^{-1}) \gamma' \beta'. $$
These relations will be used in section 3. 

Above relations are equivalent to the equation 
$$R_{p,q} T_1' T_2' = T_2' T_1' R_{p,q} \eqno(4)$$
where $T_1' = T' \otimes I$, $T_2' = I \otimes T'$ and 
$$R_{p,q} 
 = \left(\matrix{ q & 0          & 0 & 0 \cr 
                  0 & qp^{-1}    & 0 & 0 \cr
                  0 & q - p^{-1} & 1 & 0 \cr 
                  0 & 0          & 0 & p^{-1}  \cr }
\right). \eqno(5)$$
Here we employ the convenient grading notation 
$$(T_1)^{ij}{}_{kl} 
   = (T \otimes I)^{ij}{}_{kl} 
   = (-1)^{k(j+l)} T^i{}_k \delta^j{}_l, \eqno(6) $$
$$ (T_2)^{ij}{}_{kl} 
   = (I \otimes T)^{ij}{}_{kl} 
   = (-1)^{i(j+l)} T^j{}_l \delta^i{}_k. \eqno(7)$$

\noindent
{\bf 3. The two-parametric deformation of GL$(1|1)$} 

\noindent
We introduce new coordinates $x$ and $\xi$ by 
$$U = g^{-1}_{_{h_1}} U', \quad 
  U' = \left(\matrix{x' \cr \xi'\cr}\right) \eqno(8) $$
where 
$$ g_{_{h_1}} = \left(\matrix{1 & 0 \cr f_1 & 1 \cr} \right), \quad 
    f_1 = {{h_1}\over {p - 1}} \eqno(9) $$
Here the deformation parameter $h_1$, in contrast to the usual 
situation, is an odd (grassmann) number which has the following 
properties 
$$h_1^2 = 0 ~~\mbox{and}~~ h_1 \xi = - \xi h_1. \eqno(10)$$
Now, in the limit $p \rightarrow 1$ we get the following exchange 
relations 
$$x \xi = \xi x + h_1 x^2, \quad \xi^2 = - h_1 x \xi. \eqno(11)$$
These relations define a new deformation, which we called the 
$h_1$-deformation, of the algebra of functions on the Manin 
superplane generated by $x$ and $\xi$, and we denote it by $A_{h_1}$. 

Let us now consider other (dual) coordinates $\eta$ and $y$ with 
$$V = g^{-1}_{h_2} V', \quad 
  V' = \left(\matrix{\eta' \cr y'\cr}\right) \eqno(12) $$
where 
$$ g_{_{h_2}} = \left(\matrix{1 & f_2 \cr 0 & 1 \cr} \right), \quad 
    f_2 = {{h_2}\over {q - 1}} \eqno(13) $$
Again, the deformation parameter $h_2$ is an odd (grassmann) number 
and it has the following properties 
$$h_2^2 = 0 ~~\mbox{and}~~ h_2 \eta = - \eta h_2. \eqno(14)$$
Next, taking the $q \rightarrow 1$ limit we obtain the following 
relations, which define the dual $h_2$-superplane $A^*_{h_2}$ as generated 
by $\eta$ and $y$ with the commutation rules 
$$\eta^2 = - h_2 \eta y, \quad \eta y = y \eta - h_2 y^2. \eqno(15)$$
Note that in order to obtain the superplane $A_{h_1}$ and its dual 
$A^*_{h_2}$, we introduced above, two matrices $g_{h_1}$ and $g_{h_2}$. 
Of course, this result could be obtained by using a single matrix 
$g = g_{h_1} g_{h_2}$. But in that case the required steps are rather 
complicated and tedious. However, in section 4, in order to obtain an 
R-matrix we shall take $g = g_{h_1} g_{h_2}$. 

We now consider the linear transformations 
$$ T : A_{h_1} \longrightarrow A_{h_1} ~~\mbox{and}~~ 
   T : A_{h_2}^* \longrightarrow A^*_{h_2}. \eqno(16)$$
Then, we define the corresponding $(h_1,h_2)$-deformation of the 
supergroup GL$(1|1)$ as a quantum matrix supergroup GL$_{h_1,h_2}(1|1)$ 
generated by $a$, $\beta$, $\gamma$, $d$ which satisfy the following 
$(h_1,h_2)$-commutation relations 
$$a \beta = \beta a - h_2 (a^2 - \beta \gamma - ad), \quad 
  d \beta = \beta d + h_2 (d^2 + \beta \gamma - da), $$
$$a \gamma = \gamma a + h_1 (a^2 + \gamma \beta - ad), \quad 
  d \gamma = \gamma d - h_1 (d^2 - \gamma \beta - da), $$
$$ \beta^2 = h_2 \beta (a - d), \quad \gamma^2 = h_1 \gamma (d - a), $$
$$ \beta \gamma = - \gamma \beta + (h_1 \beta - h_2 \gamma)(d - a), \eqno(17)$$
$$ ad = da + (h_1 \beta + h_2 \gamma)(a - d) - h_1 h_2 (a^2 - 2 da + d^2) $$
provided that $\beta$ and $\gamma$ anticommute with $\xi$, $\eta$, $h_1$ and 
$h_2$, and 
$$h_1 h_2 = - h_2 h_1. \eqno(18)$$
One can see that when $h_2 = 0$, these relations go back to those of 
GL$_h(1|1)$ in Ref. 13. 

Alternativelly, the relations (17) can be obtained using the following 
similarity transformation which was used first in [12]:  
$$T' = g T g^{-1} \eqno(19)$$
where in our case 
$$g = g_{h_1} g_{h_2}. \eqno(20)$$
To do this, we use the relations (3) and then take the limits 
$p \rightarrow 1$, $q \rightarrow 1$. 

We denote by ${\cal A}_{h_1,h_2}$ the algebra generated by the elements $a$, 
$\beta$, $\gamma$, $d$ with the relations (17). The algebra 
${\cal A}_{h_1,h_2}$ is a (graded) Hopf algebra with the usual co-product 
$$\Delta(t^i{}_j) = t^i{}_k \otimes t^k{}_j \eqno(21)$$
(sum over repeated indices), co-unit 
$$\varepsilon(t^i{}_j) = \delta^i{}_j, \eqno(22)$$
and the antipode (co-inverse), which is the same as in [13], 
$$T^{-1}  
 = \left(\matrix{   a^{-1} + a^{-1} \beta d^{-1} \gamma a^{-1} & 
                  - a^{-1} \beta d^{-1} \cr \cr
                  - d^{-1} \gamma a^{-1} & 
                    d^{-1} + d^{-1} \gamma a^{-1} \beta d^{-1} \cr}
\right), \eqno(23) $$
provided that the formal inverses $a^{-1}$ and $d^{-1}$ exist. 

The quantum superdeterminant of $T$ is defined as, like that in the 
quantum supergroup GL$_h(1|1)$, 
$$D_{h_1,h_2} = ad^{-1} - \beta d^{-1} \gamma d^{-1} \eqno(24)$$
which is independent of the relations (17). The equation 
$$ad^{-1} - \beta d^{-1} \gamma d^{-1} = 
d^{-1}a - d^{-1} \beta d^{-1} \gamma \eqno(25)$$
is also valid, however the proof is rather lengthly but straightforward. 
It can be checked using the relations 
$$d^{-1} \beta = \beta d^{-1} - h_2 (1 - ad^{-1} + 
  d^{-1} \beta \gamma d^{-1}),$$
$$d^{-1} \gamma = \gamma d^{-1} + h_1 (1 - ad^{-1} - 
  d^{-1} \gamma \beta d^{-1}),$$
$$a d^{-1} = d^{-1} a + h_1 d^{-1} \beta (1 - ad^{-1}) + 
  h_2 (1 - d^{-1}a) \gamma d^{-1},$$
$$\gamma d^{-1} \gamma = 0, \quad 
  h_1 \beta d^{-1} \gamma \beta = - h_1 h_2 \beta \gamma (ad^{-1} -1). $$
It can be also verified that $D_{h_1,h_2}$ commutes with all matrix elements 
of $T$ (and $h_1$, $h_2$), that is, $D_{h_1,h_2}$ 
belongs to the centre of the algebra 
$$T D = DT.$$
Moreover, it can be checked that $D_{h_1,h_2}$ has the multiplicative property 
$$\Delta(D_{h_1,h_2}) = D_{h_1,h_2} \otimes D_{h_1,h_2}. \eqno(26)$$

We close this section with the following note. If we set $h_1 h_2 = 0$ 
then only the last term in the last relation in (17) vanishes. 
Essentially, we can eliminate the factor $h_1 h_2$ from the last 
equation in (17). Indeed some algebra gives 
$$ad = da + h_1 \beta (a - d) + h_2 (a - d) \gamma. \eqno(27) $$
Therefore the factor $h_1 h_2$ does not appear in any of the relations 
in (17). Thus we can demand that $h_1 h_2 = 0$. In this situation, 
the relations (17) can be easily obtained from (19). This issue will be 
used in the following two sections. 

\noindent
{\bf 4. R-matrix for GL$_{h_1,h_2}(1|1)$}

\noindent
We shall obtain an R-matrix for the quantum supergroup GL$_{h_1,h_2}(1|1)$ 
from the R-matrix of GL$_{p,q}(1|1)$. We know that the associative algebra 
(3) is equivalent to 
$$R_{p,q} T'_1 T_2' = T_2' T_1' R_{p,q} $$
[see equ.s (4)-(7)]. Now substituting (19) into (4) and defining the R-matrix 
$R_{h_1,h_2}$ as 
$$R_{h_1,h_2} = \lim_{p \rightarrow 1} \lim_{q \rightarrow 1} 
  (g \otimes g)^{-1} R_{p,q}(g \otimes g), \eqno(28)$$
where the matrix $g$ is given by (20), 
we get the following R-matrix $R_{h_1,h_2}$ 
$$R_{h_1,h_2}  
 = \left(\matrix{ 1 - h_1h_2 & - h_2     & h_2        & 0 \cr 
                  - h_1      &   1       & - h_1 h_2  & h_2 \cr
                    h_1      & - h_1 h_2 &   1        & h_2 \cr 
                    0        &   h_1     &   h_1      & 1 + h_1 h_2  \cr }
\right) \eqno(29)$$
which gives the $(h_1,h_2)$-deformed algebra of functions on 
GL$_{h_1,h_2}(1|1)$ with the equation 
$$R_{h_1,h_2} T_1 T_2 = T_2 T_1 R_{h_1,h_2}. \eqno(30)$$

Note that the above R-matrix for the special case $h_2 = 0$ coincides 
with the one in Ref. 13. For the special case $h_1 = 0 = h_2$, the 
above R-matrix becomes the unit matrix. Also 
$\widehat{R}_{h_1,h_2} = P R_{h_1,h_2}$, where $P$ is 
the super permutation matrix, satisfies 
$$\widehat{R}_{h_1,h_2}^2 = I,$$
and thus it has two eigenvalues $\pm 1$. 

If we set $h_1h_2 = 0$ in (29) then the matrix $R_{h_1,h_2}$ in (29) can be 
decomposed in the form 
$$R_{h_1,h_2} = R_{h_1} R_{h_2} \eqno(31)$$
where 
$$R_{h_1}  
 = \left(\matrix{   1   &   0   & 0   & 0 \cr 
                  - h_1 &   1   & 0   & 0 \cr
                    h_1 &   0   & 1   & 0 \cr 
                    0   &   h_1 & h_1 & 1  \cr }\right), \quad 
R_{h_2}  
 = \left(\matrix{ 1  & - h_2 & h_2 & 0 \cr 
                  0  &   1   & 0   & h_2 \cr
                  0  &   0   & 1   & h_2 \cr 
                  0  &   0   & 0   & 1   \cr }\right). \eqno(32)$$
Here the matrix $R_{h_1}$ coincides with the $R_{h}$ matrix of [13]. 

It can be checked that these matrices both satisfy the graded 
Yang-Baxter equation 
$$R_{12} R_{13} R_{23} = R_{23} R_{13} R_{12}, \quad 
  R \in \{R_{h_1}, R_{h_2}\} \eqno(33)$$
where 
$$(R_{12})^{abc}{}_{def} = R^{ab}{}_{de} \delta^c{}_f, $$
$$(R_{13})^{abc}{}_{def} = (-1)^{b(c + f)}R^{ac}{}_{df} \delta^b{}_e, $$
$$(R_{23})^{abc}{}_{def} = (-1)^{a(b + c + e + f)}R^{bc}{}_{ef} \delta^a{}_d. 
  \eqno(34)$$
Also, the matrix $R_{h_2}$ obeys the ungraded Yang-Baxter equation. 

If we set $\widehat{R} = P R$ then one can show that they satisfy 
the graded braid equation 
$$\widehat{R}_{12} \widehat{R}_{23} \widehat{R}_{12} = 
  \widehat{R}_{23} \widehat{R}_{12} \widehat{R}_{23}, \quad 
  \widehat{R} \in \{\widehat{R}_{h_1}, \widehat{R}_{h_2}\} \eqno(35)$$
where $\widehat{R}_{h_1} = P R_{h_1}$ and $\widehat{R}_{h_2} = P R_{h_2}$ 
with the grading again given by (34). Note that the matrix $R_{h_2}$ 
does not satisfy the ungraded braid equation. 

\vfill\eject
\noindent
{\bf 5. Differential Calculus on the Quantum Superplane} 

\noindent
It is well known, through the work of Woronowicz [16], that quantum 
groups provide a concrete example of non-commutative differential 
geometry. Wess and Zumino [17] developed a differential calculus on 
the quantum hyperplane covariant with respect to quantum groups. They 
have shown that one can define a consistent differential calculus on 
the non-commutative space of the quantum hyperplane. 

In this section we shall construct a covariant differential calculus on 
the quantum superplane. Before disscussing the differential calculus, we 
note the following. Let us consider the following dual (exterior) 
superplane $\Lambda_q$ as generated by $\varphi'$ and $u'$ with the 
relations 
$$\varphi'^2 = 0, \quad \varphi' u' + q^{-1} u' \varphi' = 0. \eqno(36)$$
We define 
$$\widehat{V} = g^{-1}_{_{h_2}} \widehat{V}' \eqno(37)$$
where the matrix $g_{_{h_2}}$ is given by (13). Then we get the following 
relations 
$$\varphi^2 = h_2 \varphi u, \quad u \varphi + \varphi u = - h_2 u^2 \eqno(38)$$
under the assumption 
$$h_2 u = - u h_2. \eqno(39) $$
These relations define a deformation of the (exterior) algebra of 
functions on the dual (exterior) superplane generated by $\varphi$ and 
$u$, and we denote it by $\Lambda_{h_2}$. 

Now one can check that the transformations 
$$T : A_{h_1} \longrightarrow A_{h_1}, \quad 
  T : \Lambda_{h_2} \longrightarrow \Lambda_{h_2} \eqno(40)$$
define a two-parameter deformation of the supergroup GL$(1|1)$, that is, 
they give the $(h_1,h_2)$- commutation relations in (17). 

We now pass to the differential calculus on the superplane. Consider the 
coordinates $x$ and $\xi$, belonging to the associative algebra 
$A_{h_1}$ which satisfy the commutation relations 
$$U^i U^j = \left(\widehat{R}_{h_1}\right)^{ij}{}_{kl} U^k U^l, \quad 
 \widehat{R}_{h_1} = P R_{h_1} \eqno(41)$$
where $U = (x, \xi)^T$. These relations are equivalent to (11). Similarly, 
one can express the relations (38) in the form 
$$\widehat{V}^i \widehat{V}^j = - \left(\widehat{R}_{h_2}\right)^{ij}{}_{kl} 
  \widehat{V}^k \widehat{V}^l. \eqno(42)$$ 
Note that the relations (41) and (42) can be also expressed using the 
$R_{h_1,h_2}$ matrix in (29) provided it is permuted by the super 
permutation matrix $P$. 

Denoting the partial derivatives with 
$$\partial_i = {\partial\over {\partial U^i}},   \quad 
  \widehat{\partial}_i = {{\widehat{\partial}}\over {\widehat{\partial} 
  \widehat{V}^i}},  \eqno(43)$$
one arrives at 
$$\partial_j U^i = \delta^i{}_j + 
   \left(\widehat{R}_{h_1,h_2}\right)^{ik}{}_{jl} U^l \partial_k \eqno(44)$$
and 
$$\widehat{\partial}_j \widehat{V}^i =  
   \left(\widehat{R}_{h_1,h_2}\right)^{ik}{}_{jl} \widehat{V}^l 
   \widehat{\partial}_k. \eqno(45)$$
Then we set up the relations between the coordinates of $U$, $\widehat{V}$ 
and their partial derivatives as follows: [Here, for simplicity we assumed 
that $h_1 h_2 = 0$ ]
$$\partial_x x = 1 + x \partial_x + h_1 x \partial_\xi - h_2 \xi \partial_x, 
  \quad 
 \partial_x \xi = \xi \partial_x - h_1 (x \partial_x + \xi \partial_\xi), $$
$$\partial_\xi x = x \partial_\xi + h_2 (x \partial_x + \xi \partial_\xi), 
  \quad 
  \partial_\xi \xi = 1 - \xi \partial_\xi - h_1 x \partial_\xi + 
  h_2 \xi \partial_x, \eqno(46)$$
and from (45) 
$$\partial_\varphi \varphi = \varphi \partial_\varphi + 
  h_1 \varphi \partial_u - h_2 u \partial_\varphi, \quad 
  \partial_\varphi u = u \partial_\varphi - h_1 (\varphi \partial_\varphi + 
  u \partial_u), $$
$$\partial_u \varphi = \varphi \partial_u + h_2 (\varphi \partial_\varphi + 
  u \partial_u), \quad 
  \partial_u u = - u \partial_u - h_1 \varphi \partial_u + 
  h_2 u \partial_\varphi. \eqno(47)$$

The relations between the coordinates $U$ and $\widehat{V}$ are 
$$U^i \widehat{V}^j = \left(\widehat{R}_{h_1,h_2}\right)^{ij}{}_{kl} 
  \widehat{V}^k U^l, \eqno(48)$$ 
which read 
$$x\varphi = \varphi x + h_2 (ux - \varphi \xi), \quad 
  xu = ux + h_1 \varphi x + h_2 u \xi$$
$$\xi \varphi = \varphi \xi - h_1 \varphi x + h_2 u \xi, \quad 
 \xi u = - u \xi - h_1 (\varphi \xi + ux). \eqno(49)$$

Finally one gets from the following equation 
$$\partial_i \partial_j = 
  \left(\widehat{R}_{h_1,h_2}\right)^{lk}{}_{ji} 
  \partial_k \partial_l \eqno(50)$$
which are the commutation relations among the partial derivatives: 
$$\partial_x \partial_\xi = \partial_\xi \partial_x - h_2 \partial_x^2, $$
$$\partial_\xi^2 = h_2 \partial_\xi \partial_x. \eqno(51)$$

One of the interesting problems is to construct $U_{h_1,h_2}(gl(1|1))$. 
Work on this issue is in progress. 

\noindent
{\bf Acknowledgement}

\noindent
This work was supported in part by {\bf T. B. T. A. K.} the 
Turkish Scientific and Technical Research Council. 

\def\refname{References}


\baselineskip=12pt
\begin{thebibliography}{99}
\bibitem{D:gnus} Drinfeld, V. G., {\em Quantum groups}, 
    in {\it Proc. } IMS, Berkeley, 1986.
\bibitem{RTF:gnus} Reshetikhin, N. Y, Takhtajan, L. A. and Faddeev, L. D., 
    {\it Leningrad Math. J. } {\bf 1} (1990), 193-225.
\bibitem{M:gnus} Manin, Yu I., 
    {\it Commun. Math. Phys.} {\bf 123} (1989), 163-175. 
\bibitem{DMMZ:gnus} Demidov, E.E, Manin, Yu I., Mukhin, E. E and 
    Zhdanovich, D.V., 
    {\it Prog. Theor. Phys. Suppl.} {\bf 102} (1990), 203-218.
\bibitem{EOW:gnus} Ewen, H., Ogievetsky, O. and Wess, J., 
    {\it Lett. Math. Phys.} {\bf 22} (1991), 297-305. 
\bibitem{Z:gnus} Zakrzewski, S., 
    {\it Lett. Math. Phys.} {\bf 22} (1991), 287-289. 
\bibitem{W:gnus} Woronowicz, S. L., 
    {\it Rep. Math. Phys.} {\bf 30} (1991), 259-269. 
\bibitem{O:gnus} Ohn, C. H., 
    {\it Lett. Math. Phys.} {\bf 25} (1992), 85-88. 
\bibitem{Ku:gnus} Kupershmidt, B. A., 
    {\it J. Phys. A: Math. Gen.} {\bf 25} (1992), L1239-L1244. 
\bibitem{Ka:gnus} Karimipour, V., 
    {\it Lett. Math. Phys.} {\bf 30} (1994), 87-98. 
\bibitem{A:gnus} Aghamohammadi, A., 
    {\it Mod. Phys. Lett. Math. A} {\bf 8} (1993), 2607-2613.
\bibitem{AKS:gnus} Aghamohammadi, A., Khorrami, M. and Shariati, A. 
    {\it J. Phys. A: Math. Gen.} {\bf 28} (1995), L225-231.
\bibitem{DP:gnus} Dabrowski, L. and Parashar, P., 
    to appear in {\it Lett. Math. Phys.} 
\bibitem{DP:gnus} Dabrowski, L. and Wang, L. {\it Phys. Lett. B} {\bf 266} 
    (1991), 51-54. 
\bibitem{CC:gnus} Chakrabarti, R. and Jagannathan, R., {\it J. Phys. A: Math. 
    Gen.} {\bf 24} (1991), 5683-5701. 
\bibitem{SL:gnus} Woronowicz, S. L., {\it Commun. Math. Phys.} {\bf 122} 
    (1989), 125-170. 
\bibitem{WZ:gnus} Wess, J. and Zumino, B., {\it Nucl. Phys. B } (Proc. Suppl.) 
    {\bf 18} B (1990), 302.

\end{thebibliography}
\end{document}